\documentclass[a4paper,10pt,english]{smfart}
\usepackage{amsfonts}
\usepackage{amsmath} 
\usepackage{latexsym}
\usepackage{array}
\usepackage{amssymb}
\usepackage{smfthm}
\usepackage{bull}
\theoremstyle{plain}

\newcommand{\R}{  \mathbb{R}   }
\newcommand{\eps}{\varepsilon}
\newcommand{\e}{  \text{e}   }

\newcommand{\N}{  \mathbb{N}   }

\author{Laurent Thomann}
\address{Universit\'e Paris-Sud, Math\'ematiques, B\^at 425\\Tel 0169155785\\ 91405
Orsay Cedex.}
\email{ laurent.thomann@math.u-psud.fr}
\urladdr{http://www.math.u-psud.fr/~thomann}
\title[ Instability for the Gross-Pitaevski equation]{ Geometric and projective instability for the Gross-Pitaevski equation} 
 \date{}

\begin{document}
\frontmatter
\begin{abstract}
Using variational methods, we construct approximate solutions for the Gross-Pitaevski equation 
which concentrate on circles in $\R^3$. These solutions will help to show that the $L^2$ flow is unstable for the 
usual topology and for the projective distance.
\end{abstract}

\maketitle
\mainmatter
\section{Introduction}
In this paper we deal with the equations 
\begin{equation}\label{grosspita}
\left\{
\begin{aligned}
&i h\partial_t u+h^2\Delta u-|x|^2u   = a_hh^2|u|^{2} u, \quad (t,x)\in \R^{1+3},\\
&u(0,x)= u_0(x) \in L^2(\R^3),
\end{aligned}
\right.
\end{equation}
where $h>0$ is a small parameter and $a_h$ a constant which depends on
$h$, that can be either positive (defocusing case) 
or negative (focusing case). In all the paper we assume that there
exists a constant $A>0$, independent of $h$, such that $|a_h|\leq
A$.\\
 This equation appears in the study of Bose-Einstein condensates; for more details 
see \cite{Pitaevskii}.\\
In the following we will refer to the definitions:
\begin{defi}
(Geometric instability) We say that the Cauchy problem $(\ref{grosspita})$ is geometrically unstable if 
there exist $u_h^1, u_h^2\in L^2(\R^3)$ solutions of $(\ref{grosspita})$ with initial data 
$u_h^1(0), u_h^2(0)\in L^2(\R^3)$ such that $\|u_h^1(0)\|_{L^2},
\|u_h^2(0) \|_{L^2} \leq C$ where $C$ is a constant independent of
  $h$ and $t_h>0$ such that 
$$\frac{ \|(u_h^2-u_h^1)(t_h) \|_{L^2}}{ \|(u_h^2-u_h^1)(0) \|_{L^2}} \longrightarrow +\infty \;\; 
\text{when} \;\;h\longrightarrow 0.$$
\end{defi}

\begin{defi}
(Projective instability) We say that the Cauchy problem $(\ref{grosspita})$ is projectively unstable if 
there exist $u_h^1, u_h^2\in L^2(\R^3)$ solutions of $(\ref{grosspita})$ with initial data 
$u_h^1(0), u_h^2(0)\in L^2(\R^3)$  such that $\|u_h^1(0)\|_{L^2},
\|u_h^2(0) \|_{L^2} \leq C$ where $C$ is a constant independent of
  $h$ and $t_h > 0$ such that 
$$\frac{ d_{\text{pr}}\left(u_h^2(t_h),u_h^1(t_h)\right)}{d_{\text{pr}}\left( u_h^2(0),u_h^1(0) \right)} \longrightarrow +\infty \;\; 
\text{when} \;\;h\longrightarrow 0.$$
Here $ d_{\text{pr}}$ denotes the complex projective distance defined by 
$$ d_{\text{pr}}(v_1,v_2)=\arccos{\left(\frac{\left|\left<v_1,v_2 \right>\right|}{\|v_1\|_{L^2}\|v_2\|_{L^2}}\right)}
\;\;for   \;\; v_1,v_2 \in L^2(\R^3). $$
\end{defi}

\begin{enonce}{Notations}
In this paper $c$, $C$ denote constants the value of which may change
from line to line. These constants will always be independent of $h$. We use the notations $a\sim b$,  
$a\lesssim b$, $a\gtrsim b$, if $\frac1C b\leq a\leq Cb$ , $a\leq Cb$, $b\leq   Ca$
respectively. We write $a\ll b$, $a\gg b$ if $a\leq Kb$, $a\geq K b$ for some large constant
$K$ which is  independent of $h$.
\end{enonce}

The first result of this paper is

\begin{theo}\label{thmgeom}
Let $h^{-1}\in \N$. In each of the following cases, there exist
 $c_0>0$ and $u_h^1, u_h^2\in L^2(\R^3)$ solutions
 of $(\ref{grosspita})$ with initial data 
 $ \|u_h^2(0)\|_{L^2}$, $ \|u_h^1(0) \|_{L^2}\to \kappa$ such that if
 $|a_h|\kappa^2\leq c_0$, we have:\\
$(i)$ Assume $a$ is independent of $h$ and $\kappa |a| t \gg 1$,
$$\frac{ \|(u_h^2-u_h^1)(t) \|_{L^2}}{ \|(u_h^2-u_h^1)(0) \|_{L^2}} \gtrsim |a| \kappa t.$$
$(ii)$  Assume $ |a_h| t_h \longrightarrow +\infty$ when
 $h\longrightarrow 0$ with $t_h\ll\log{\frac1h}$, then
$$\sup_{0\leq t \leq t_h}\|(u_{h}^2-u_{h}^1)(t) \|_{L^2} \gtrsim 1, $$
 but
$$\|(u_{h}^2-u_{h}^1)(0) \|_{L^2} \longrightarrow 0.$$ 
In particular, the Cauchy problem $(\ref{grosspita})$
is geometrically  unstable
\end{theo}

\noindent Denote by $x=(x_1,x_2,x_3)$ the current point in $\R^3$. In cylindrical coordinates
$(x_1=r\cos\theta,x_2=r\sin\theta,x_3=y)$, the
functions considered in Theorem \ref{thmgeom} take the form
\begin{equation}\label{data}
u_h(0,x)= \kappa_h h^{-\frac12}\text{e}^{i\frac{k^2}{h}\theta}
{v_0}\big(\frac{r-k}{\sqrt{h}},\frac{{y}}{\sqrt{h}}\big),
\end{equation}
where $k\in\N$, $v_0\in L^2(\R^2)$ and
\begin{equation}\label{sol}
u(t,x)=u_h(0,x)\e^{-i\lambda_h t}+w_h(t,x),\end{equation}
with $w_h$ a small error term in $ L^2(\R^3)$, at least for times when
instability effects occur.\\
The Ansatz \eqref{data} shows that the function $u$ in \eqref{sol} will
concentrate on the circle $(x_1^2+x_2^2=k^2, x_3=0)$ in $\mathbb{R}^3$.\\
To prove Theorem \ref{thmgeom}, we consider two initial data of the
form \eqref{data} associate with $\kappa$ and $\kappa'$ such that
$|\kappa'-\kappa|$ is small, and therefore the initial data are close
in the $L^2$-norm, but we will see that the solutions do not remain
close to each other after a time $t$. \\
The construction of two solutions  to \eqref{grosspita} of the form
\eqref{data},\eqref{sol} which
concentrate on disjoint circles yield the following result 
\begin{theo}\label{thmproj}
Let $h^{-1}\in \N$. There exist $c_0>0$ and $u_h^1, u_h^2\in L^2(\R^3)$ solutions
 of $(\ref{grosspita})$ with initial data 
 $ \|u_h^2(0)\|_{L^2}$, $ \|u_h^1(0) \|_{L^2}\to \kappa$ such that if
 $|a_h|\kappa^2\leq c_0$ and $ |a_h|t_h \longrightarrow
 +\infty$ when $h\longrightarrow 0$ with $t_h\ll\log{\frac1h}$, we have
$$\sup_{0\leq t \leq
 t_h}d_{\text{pr}}\left(u_{h}^2(t),u_{h}^1(t)\right)  \gtrsim 1, $$
 but
$$d_{\text{pr}}\left(u_{h}^2(0),u_{h}^1(0)\right) \longrightarrow 0.$$ 
In particular, the Cauchy problem $(\ref{grosspita})$
is  projectively unstable.
\end{theo}
\noindent The part $(i)$ of Theorem \ref{thmgeom} shows that there is no Lipschitz dependence between the solutions of equation 
\eqref{grosspita} and the initial data in the regime $\kappa a t \gg 1$, whereas the part $(ii)$ and Theorem \ref{thmproj} assert 
that the dependence is not continuous, but for larger times. Both types of instabilities are nonlinear  behaviour, but the first 
one is weaker than the second.\\[5pt]
The instability results of Theorem \ref{thmgeom} are not new in the case
$a>0$. R. Carles \cite{Carles1} shows the instability, for finite
times, of the equation 
$$i h\partial_t v+h^2\Delta v-|x|^2v   = f(h^k|v|^{2}) v,\quad  (t,x)\in\R^{1+n}, $$
when $n\geq 2$, $1<k<n$, and $f\in\mathcal{C}^{\infty}(\R_+,\R)$ with $f'>0$.
\\[5pt]
In \cite{BGT1}, N. Burq, P. G\'erard and N. Tzvetkov have  pointed out geometric instability for the cubic 
Schr\"odinger equation $i \partial_t u+\Delta_{\mathbb{S}^2} u  = a|u|^{2}u $
 on $\mathbb{S}^2$ when $a>0$. This phenomenon doesn't occur on $L^2(\R^3)$ for the equation 
$i \partial_t u+\Delta u   = a|u|^{2}u $ in $L^2(\R^3) $, it is therefore strongly related to the geometry of
the operator and of the manifold we work on. Here there is no
semiclassic parameter in the equations, but we could obtain similar
results in this  latter case with a scaling argument, as these instability
effects are local. There are stronger instability phenomenona in $H^s$
norm, for $0<s<\frac12$ or for $s$ negative, for more details see
\cite{CCT2} or \cite{CCT1} for the one dimensional case.\\[5pt]
In \cite{BZ}, N. Burq and M. Zworski prove Theorem $\ref{thmgeom}$ in the case $a>0$. To obtain geometric instability, 
they expand the solution on the Hilbertian basis given by the eigenfunctions of $-h^2\Delta +|x|^2$. The nonlinear term 
in  $(\ref{grosspita})$ induces a phase shift in time for the groundstate and this yields the result. We will give 
a more precise description of the solution by solving a pertubated eigenvalue problem for the harmonic oscillator and 
this will also treat the focusing case.
They also obtain projective instability for the equation
\begin{equation*}
i h\partial_t u+h^2\Delta u-V(x)u   = ah^2|u|^{2} u,
\end{equation*}
where $V$ is a cylindrically symmetric potential with respect to the variable $y=x_3$, but they have to add 
the following assumption: Denote by $r=\sqrt{x_1^2+x_2^2}$ then the function $(r,y)\longmapsto V(r,y)+r^{-2}$ 
has two distinct absolute non-degenerate minima $(r_j,y_j), j=1,2, $ and its Hessian at  $(r_j,y_j)$ are 
equal. We use  a variational method to construct quasimodes which are localized  on circles in $\R^3$, which allows to remove such an hypothesis. This idea comes from an unpublished work from  N. Burq, P. G\'erard and N. Tzvetkov .
\\\\
\noindent Thanks to the form $F(|u|^2)u$ of the nonlinearity in
\eqref{grosspita}, we look for a solution $u$ which writes
$u(t,x)=\e^{-i\lambda t}f(x)$. Then $f$ has to satisfy
$$\left(-h^2\Delta+|x|^2\right)f=h\lambda f -a_hh^2|f|^2f.$$
In the case $a_h=0$, $f$ is an eigenvector of the operator
$-h^2\Delta+|x|^2$ associate with the eigenvalue $h\lambda$. In the
general case, the term 
$a_hh^2|f|^2f$ will be treated as a perturbation of the linear problem
$$\left(-h^2\Delta+|x^2|\right)f=h\lambda f.$$
In fact, we will find a development in powers of $h$ of $h\lambda$
and $f$
$$h\lambda\sim\sum_{k\geq0}\mu_kh^k,\quad f\sim\sum_{k\geq0}f_kh^k, $$
by solving a cascade of equations. This will be done in cylindrical
coordinates: Write $x=(x_1,x_2,x_3)$ and make the cylindrical change of variables $x_1  = r \cos \theta$, 
$x_2  = r \sin \theta$ and $x_3  = y$ with $(r, \theta,y) \in \mathbb{R}_+^* \times [0,2\pi[\times \mathbb{R}.$
Then the Laplace operator takes the form

\begin{equation*}
\Delta=\frac{1}{r^2}\partial_{\theta}^2+\partial_{r}^2+\frac{1}{r}\partial_r+\partial_{y}^2.
\end{equation*}
Let $\kappa$ be a positive constant and $k$ a positive integer, we want to find a solution of $(\ref{grosspita})$ 
of the form 

\begin{equation}\label{defutilde}
\tilde{u}=\kappa h^{-\frac{1}{2}}\text{e}^{-i\lambda t}\text{e}^{i\frac{k^2}{h}\theta} \tilde{v}(r,y,h),
\end{equation}
where $\lambda$ is a constant to be determined, and $\tilde{v}$ a real function which therefore has to satisfy
\begin{equation*}
-h^2(\partial_{r}^2+\partial_{y}^2)\tilde{v}+(\frac{k^4}{r^2}+r^2+y^2)\tilde{v}=\lambda h\tilde{v} -a_hh^2\kappa^2\tilde{v}^3+h^2\frac{1}{r}\partial_r\tilde{v}.
\end{equation*}
Notice that we have to choose $h^{-1}\in \N$ so that \eqref{defutilde} makes sense for all $k\in \N$.
We try to construct $\tilde{v}$ which concentrates exponentially at the minimum of the potential 
$V=\frac{k^4}{r^2}+r^2+y^2$, i.e. at $(r,y)=(k,0)$.\\
Thus we make the change of variables $r=k+\sqrt{h}\rho$, $y=\sqrt{h}\sigma$
and set $ \tilde{v}(r,y,h)=v(\frac{r-k}{\sqrt{h}},\frac{y}{\sqrt{h}},h).$\\
We write the Taylor expansion of $V$ in $h$:
\begin{eqnarray*}
\frac{k^4}{(k+\sqrt{h}\rho)^2}
+(k+\sqrt{h}\rho)^2+h\sigma^2 &=&
2k^2+(4\rho^2+\sigma^2)h-\frac4k\rho^3h^{\frac{3}{2}} \\
& & +\frac{5}{k^2}\rho^4h^2 +R(\rho,h)h^{\frac{5}{2}}.
\end{eqnarray*}
Then $v$ has to be solution of 

\begin{eqnarray}
Eq(v)&:=&P_0v-
\frac{\lambda h-2k^2}{h}v +a_h\kappa^2{v}^3-
h^{\frac{1}{2}}(\frac{1}{k+\sqrt{h}\rho}\partial_{\rho}v+\frac4k\rho^3v) \nonumber \\
 & & +\frac{5}{k^2}\rho^4hv - h^{\frac{3}{2}}Rv=0\label{eqv},
\end{eqnarray}
where $P_0=-(\partial_{\rho}^2+\partial_{\sigma}^2)+(4\rho^2+\sigma^2).$
 Now, write

\begin{eqnarray*}
v(\rho,\sigma,h)&=&v_0(\rho,\sigma)+h^{\frac{1}{2}}v_1(\rho,\sigma)+hv_2(\rho,\sigma)
+h^{\frac{3}{2}}w(\rho,\sigma,h) \\
\frac{\lambda h-2k^2}{h}&=&E_0+h^{\frac{1}{2}}E_1+hE_2+h^{\frac32}E_3(h).
\end{eqnarray*}
 By identifying the powers of $h$ we obtain the following equations:

\begin{eqnarray}
P_0v_0&=& E_0 v_0 -a_h\kappa^2{v_0}^3, \label{eqv_0} \\
P_0v_1&=& E_0 v_1 +E_1v_0-3a_h\kappa^2{v_0}^2v_1
+\frac1k \partial_{\rho} v_0+\frac4k\rho^3v_0,\label{eqv_1} \\
P_0v_2&=& E_0 v_2 +E_1v_1+E_2v_0-3a_h\kappa^2({v_0}^2v_2+v_0v_1^2)+\frac1k \partial_{\rho} v_1 \nonumber \\
 & & +\frac4k\rho^3v_1-\frac{1}{k^2}\rho \partial_{\rho}v_0  -\frac{5}{k^2}\rho^4v_0. \label{eqv_2} 
\end{eqnarray}

\begin{rema}
In the sequel we only mention the dependence in $k$, $\kappa$ and $a$
of the $v_j$ and $E_j$ when necessary. Moreover we write $a=a_h$.
\end{rema}

\section{Construction of the quasimodes}

\begin{prop}\label{propv_0}
There exists a constant $c_0>0$ such that if $|a|\kappa^2\leq c_0$, 
there exist $E_0>0$ and
$v_0 \in L^2(\R^2)$ satisfying $v_0 \geq0$ and $\|v_0 \|_{L^2(\R^2)}=1$, which solve $(\ref{eqv_0})$.
\end{prop}
\noindent For $\psi\in\mathcal{S}'(\R^2)$, denote by $\hat{\psi}$ its
Fourier transform, with the convention 
$$\hat{\psi}(\zeta)=\int_{\R^2}\e^{-i\zeta\cdot x}\psi(x)\text{d}x,$$
for $\psi\in L^1(\R^2)$.\\
We use a variational  method based on Rellich's criterion.

\begin{prop}
(\cite{ReedSimon4}, p 247) The set 
\begin{equation*}
S=\left\{\psi | \int_{\R^2} |\psi(x)|^2\text{d}x = 1,  \int_{\R^2}
 (1+|x|^2)|\psi(x)|^2\text{d}x \leq 1,
 \int_{\R^2} (1+|\zeta|^2)|\hat{\psi}(\zeta)|^2\text{d}\zeta \leq 1 \right\},
\end{equation*}
is a compact subset of $L^2(\R^2)$.
\end{prop}

\begin{proof}[Proof of Proposition \ref{propv_0}]
We minimize the functional
\begin{equation*}
J(u,a)=\int\left(|\nabla u|^2+(4\rho^2+\sigma^2)|u|^2+\frac{1}{2}a\kappa^2|u|^4\right),
\end{equation*}
on the space 
\begin{equation*}
H=\left\{u \in H^1(\R^2),(\rho^2+\sigma^2)^{\frac12}u\in L^2(\R^2),\|u\|_{L^2}=1 \right\}.
\end{equation*}
Now, on $H$ we have the inequality
\begin{equation*}
\|u\|_{L^4} \leq C \|u\|_{H^{\frac{1}{2}}}
 \leq C\| u\|_{ L^2}^{\frac{1}{2}}\|\nabla u\|_{ L^2}^{\frac{1}{2}}
\leq C\|\nabla u\|_{ L^2}^{\frac{1}{2}}.
\end{equation*}
Thus, there exists $c_0>0$ such that 
\begin{equation*}
\frac{1}{2}a\kappa^2\int |u|^4 \leq \frac{1}{2}\int|\nabla u|^2,
\end{equation*}
as soon as $|a|\kappa^2\leq c_0$, which we suppose from now.\\
Let $(u_n)_{n\geq 1}$ be a minimizing sequence. First, we can choose $u_n \geq 0$, because $|u_n|$ is also 
 minimizing, as $|\nabla|u_n||\leq |\nabla u_n|$. We have
\begin{equation*}
\int\left(\frac12|\nabla u_n|^2+(4\rho^2+\sigma^2)u_n^2\right) \leq J(u_n,a\kappa^2)\leq C,
\end{equation*}
with $C$ independent of $a$, $\kappa$ and $n$. We are able to apply Rellich's
criterion: there exists $v_0 \in H$ with $v_0 \geq0$ such that, up to
a subsequence, $u_n \longrightarrow v_0$, and the lower  semi-continuity 
of $J$ ensures 
\begin{equation*}
J(v_0,a\kappa^2)=\inf_{u \in H}J(u,a\kappa^2).
\end{equation*}
Then there exists a Lagrange multiplier $E_0$ such that
\begin{equation*}
P_0v_0=-(\partial_{r}^2+\partial_{y}^2){v_0}+(4\rho^2+\sigma^2){v_0}= E_0 v_0 -a\kappa^2{v_0}^3,
\end{equation*}
and $E_0$ is given by
\begin{equation*}
E_0=\int\left(|\nabla v_0|^2+(4\rho^2+\sigma^2)v_0^2+a\kappa^2{v_0}^4\right).
\end{equation*}
\end{proof}

\begin{prop}\label{propestimv_0}
Let $|a|\kappa^2\leq c_0$. There exist constants $C,c>0$ independent
of $a$, $\kappa$ such that 
for $0\leq j\leq 2$
\begin{equation}\label{estimv_0}
\left|(I-\Delta)^{\frac{j}{2}}v_0(\rho,\sigma)\right| \leq C\text{e}^{-c(|\rho|+|\sigma|)}.
\end{equation}

\end{prop}

\begin{proof}
We denote by $\xi=(\rho,\sigma)$, and we define $\varphi_{\eps}(\xi)=\text{e}^{\frac{|\xi|}{1+\eps|\xi|}}$.
 The function $\varphi_{\eps}$ is bounded and 
\begin{equation}\label{borne}
|\nabla\varphi_{\eps}| \leq\varphi_{\eps} \quad\text{a.e.}
\end{equation}
We multiply $(\ref{eqv_0})$ by $\varphi_{\eps} v_0$ and integrate over $\R^2$:
\begin{equation*}
\int \nabla(\varphi_{\eps} v_0)\nabla v_0+\int\varphi_{\eps}|\xi|^2v_0^2 \leq E_0\int \varphi_{\eps}v_0^2+
|a|\kappa^2\int \varphi_{\eps}v_0^4.
\end{equation*}
We compute $\nabla(\varphi_{\eps} v_0)= v_0
\nabla\varphi_{\eps}+\varphi_{\eps}\nabla v_0$, and use \eqref{borne}
to obtain 
\begin{equation*}
\int(\varphi_{\eps} |\nabla v_0|^2+\varphi_{\eps}|\xi|^2v_0^2) \leq E_0\int \varphi_{\eps}v_0^2+
|a|\kappa^2\int \varphi_{\eps}v_0^4+\int \varphi_{\eps}v_0|\nabla v_0|.
\end{equation*}
We set $w_0=\varphi_{\eps}^{\frac{1}{4}}v_0$, then
\begin{equation}\label{grad}
 \nabla w_0=\frac{1}{4}\varphi_{\eps}^{-\frac{3}{4}}\nabla\varphi_{\eps}v_0+
\varphi_{\eps}^{\frac{1}{4}}\nabla v_0.
\end{equation}
 From the Gagliardo-Nirenberg inequality in dimension 2
$$\|w_0\|^4_{L^4}\leq C\|w_0\|^2_{L^2}\|\nabla w_0\|^2_{L^2},$$
together with $(\ref{grad})$ we deduce
\begin{equation*}
\int \varphi_{\eps}v_0^4 \leq C \int(\varphi_{\eps}^{\frac{1}{2}}v_0^2)
\int \varphi_{\eps}^{\frac{1}{2}}(v_0^2+|\nabla v_0|^2).
\end{equation*}
As $\int v_0^2=1$ and  $\int |\nabla v_0|^2 \leq C$, Jensen's inequality gives 
\begin{eqnarray}\label{ineg1}
\int \varphi_{\eps}v_0^4& \leq & C \left(\int \varphi_{\eps}v_0^2\right)^{\frac{1}{2}}
\left(   \int \varphi_{\eps}(v_0^2+|\nabla v_0|^2)\right)^{\frac{1}{2}} \nonumber \\
   &\leq & \frac{1}{16c_0}\int \varphi_{\eps}|\nabla v_0|^2+C\int \varphi_{\eps}v_0^2.
\end{eqnarray}
We also have
\begin{equation}\label{ineg2}
\int \varphi_{\eps}v_0|\nabla v_0| \leq 
\frac{1}{4}\int \varphi_{\eps}|\nabla v_0|^2+C\int \varphi_{\eps}v_0^2.
\end{equation}
Now, write for $R>0$
\begin{equation*}
\int \varphi_{\eps}v_0^2=\int_{|\xi| <R} \varphi_{\eps}v_0^2+\int_{|\xi| \geq R} \varphi_{\eps}v_0^2
\leq \text{e}^R\int v_0^2+\frac{1}{R^2}\int \varphi_{\eps}|\xi|^2 v_0^2,
\end{equation*}
and deduce that for $R$ big enough, independent of $\eps$, there exists a constant $C$ independent of $\eps$ satisfying 
\begin{equation*}
\int \varphi_{\eps}(|\nabla v_0|^2+|\xi|^2v_0^2) \leq C.
\end{equation*}
Letting $\eps$ tend to 0  yields
\begin{equation}\label{estL2}
\e^{\frac{|\xi|}{2}}\nabla v_0 \in L^2\quad \text{and} \quad \e^{\frac{|\xi|}{2}}|\xi| v_0\in L^2.
\end{equation} 
With the help of equation $(\ref{eqv_0})$, compute
\begin{eqnarray*}
\Delta \left(v_0\e^{\frac14(\rho+\sigma)}\right)&=&a\kappa^2v_0^3\e^{\frac14(\rho+\sigma)}+
(4\rho^2+\sigma^2-E_0)v_0\e^{\frac14(\rho+\sigma)}\\
& & +\frac12(1,1)\cdot\nabla v_0\e^{\frac14(\rho+\sigma)}
+\frac{1}{16}v_0\e^{\frac14(\rho+\sigma)}.
\end{eqnarray*}
According to $(\ref{estL2})$, each term of the right hand side is in $L^2$, excepted maybe 
the first one. But denote by $\tilde{v}_0=v_0\e^{\frac{1}{12}(\rho+\sigma)}$, then $(\ref{estL2})$ shows 
that $\tilde{v}_0\in H^1(\R^2)$ and consequently $v_0 \in L^6(\R^2).$\\
Hence, with the inequality $\|w\|_{L^{\infty}}^2\lesssim \|w\|_{L^2}\|\Delta w\|_{L^2}$ applied 
to $w=v_0\e^{\frac14(\rho+\sigma)}$ we deduce $v_0\leq C\e^{-\frac14(\rho+\sigma)} $.\\
The same can be done with $\sigma$ replaced with  $-\sigma$ or $\rho$ by $-\rho$. Therefore 
$v_0\leq C\e^{-\frac14(|\rho|+|\sigma|)}$.
Equation \eqref{eqv_0} and the previous estimate give 
$$|\Delta v_0(\rho,\sigma)| \leq C\text{e}^{-c(|\rho|+|\sigma|)}.$$
To obtain the last estimation of Proposition \ref{propestimv_0}, use the interpolation inequality 
$$\|\nabla w \|^2_{L^{\infty}}\leq \| w \|_{L^{\infty}}\|\Delta w \|_{L^{\infty}},$$
applied to $w=v_0\e^{c(\pm \rho\pm \sigma)}$.
\end{proof}

\noindent We are now able to describe the behaviour of $E_0(a\kappa^2)$ and $v_0(a\kappa^2)$ when $a\kappa^2\longrightarrow 0$:

\begin{prop}\label{propconv}
\begin{equation*}
v_0(a\kappa^2)\longrightarrow
 \frac{2^{\frac{1}{4}}}{\pi^{\frac{1}{2}}}\text{e}^{-(\rho^2+\frac{1}{2}\sigma^2)}
 \quad \text{in}\;\; L^2(\R^2)\;\;\text{when}\;\;a\kappa^2 \longrightarrow 0,
\end{equation*}
and
\begin{equation}\label{dlE0}
E_0(a\kappa^2)=3+\frac{\sqrt{2}}{2\pi}a\kappa^2+o(a\kappa^2).
\end{equation}
\end{prop}

\begin{proof}
The function $u_0=\frac{2^{\frac{1}{4}}}{\pi^{\frac{1}{2}}} 
\text{e}^{-(\rho^2+\frac{1}{2}\sigma^2)}$ is 
the unique positive element in $H$ that realises the infimum of
$J(u,0)$, and is the first eigenfunction  of 
$P_0=-\Delta+(4\rho^2+\sigma^2)$ associate with the eigenvalue
$E_0(0)=3$. See \cite{Helffer2}, p 7 for details.\\
For $|a|\kappa^2\leq c_0$ we have
\begin{equation}\label{comp}
\|v_0(a\kappa^2)\|_{L^2}=1,\; \|\nabla v_0(a\kappa^2)\|_{L^2} \leq C, \; \text{and}\;\; \|\xi v_0(a\kappa^2)\|_{L^2}\leq C.
\end{equation}
By  Rellich's criterion, $(v(a\kappa^2))_{|a|\kappa^2\leq c_0}$ is
compact in $H$; let $\mathcal{A}$ be its adherence set. 
If $u \in \mathcal{A}$, there exists a sequence $b_n=a_n\kappa_n^2 \longrightarrow 0$ satisfying $v_0(b_n)\longrightarrow u$ in $L^2$. As $v_0(b_n)$ realises the infimum of $J(v,b_n)$:
\begin{equation*}
J(v_0(b_n),b_n) \leq J(u_0,b_n)=3+\frac{1}{2}b_n\int|u_0|^4,
\end{equation*}
therefore, $J(u,0)\leq 3$. As $u\geq0$, we conclude $u=u_0$, i.e. $\mathcal{A}=\{u_0\}$ and 
\begin{equation*}
v_0(a\kappa^2)\longrightarrow u_0 \quad \text{in}\;\; L^2(\R^2)\;\;\text{when}\;\;a\kappa^2 \longrightarrow 0.
\end{equation*}
Moreover $|v(a\kappa^2)|,|u_0|\leq C$, then the convergence in also in $L^4$.\\
Now, the self-adjointness of $P_0$ gives
\begin{equation*}
0=\left<(P_0-3)u_0,v(a\kappa^2)\right>=(E_0(a\kappa^2)-3)\int v(a)u_0-a\kappa^2\int v^3(a)u_0,
\end{equation*}
then from $\int v(a\kappa^2)u_0 \longrightarrow \int u_0^2=1$ and  
$\int v(a\kappa^2)^3u_0 \longrightarrow \int u_0^4=\frac{\sqrt{2}}{2\pi}$ we 
conclude $E_0(a\kappa^2)=3+\frac{\sqrt{2}}{2\pi}a\kappa^2+o(a\kappa^2).$
\end{proof}

\begin{prop}\label{propestimv_j}
Let  $|a\kappa^2|\leq c_0$. There exist $E_1,E_2 \in \R$ and
$v_1,v_2 \in L^2(\R^2)$ satisfying $v_1,v_2 \geq0$ and $\|v_1 \|_{L^2(\R^2)}$, 
$\|v_2 \|_{L^2(\R^2)}\sim 1$, which solve $(\ref{eqv_1})$ and $(\ref{eqv_2})$.\\
Moreover there exists $c>0$ such that for $l=1,2$ and $0\leq j \leq2$
\begin{equation}\label{estimv_j}
\left|(I-\Delta)^{\frac{j}{2}}v_l(\rho,\sigma)\right| \leq C\text{e}^{-c(|\rho|+|\sigma|)}.
\end{equation}
\end{prop}

\begin{proof}
Equation $(\ref{eqv_1})$ writes
\begin{equation*}
\left(P(a\kappa^2)-E_0\right)v_1=\left(-(\partial_{\rho}^2+\partial_{\sigma}^2)+V\right)v_1= E_1 v_0 +\frac1k\partial_{\rho} v_0+\frac{4}{k}\rho^3v_0,
\end{equation*}
where we denote by $P(a\kappa^2)=P_0+3a\kappa^2{v_0}^2$ and $V=4\rho^2+\sigma^2+3a\kappa^2{v_0}^2-E_0$. The potential $V$ is so that $V \longrightarrow \infty $ as $|(\rho,\sigma)|\longrightarrow \infty$, then the spectrum $\sigma(P(a))$ 
of $P(a\kappa^2)$ is purely discrete and the eigenvalues are given by the min-max principle (see \cite{ReedSimon4} p. 120).\\
The first eigenvalue of $P(a\kappa^2)$ is therefore given by 
\begin{equation*}
\mu_0(a\kappa^2)=\inf_{u\in H}\int \left( |\nabla u|^2+(4\rho^2+\sigma^2)u^2+3a\kappa^2v_0^2u^2 \right)
-E_0(a\kappa^2),
\end{equation*}
and  there exists $w_0 \in H$ with $w_0\geq0$ satisfying 
\begin{equation*}\label{eqw_0}
\left(P(a\kappa^2)-E_0\right)w_0=\left(P_0-E_0(a\kappa^2)+3a\kappa^2{v_0}^2  \right)w_0=\mu_0(a)w_0,
\end{equation*}
and one shows, as in the proof of $(\ref{propconv})$ that $w_0 \longrightarrow u_0$ in $L^2\cap L^4$.\\
Multiply $(\ref{eqv_0})$ by $u_0$ and integrate
\begin{equation*}
3a\kappa^2\int {v_0}^2w_0u_0 +(3-E_0(a\kappa^2))\int w_0u_0=\mu_0(a\kappa^2) \int w_0u_0,
\end{equation*}
then according to  $(\ref{dlE0})$, $\mu_0(a\kappa^2)\sim \frac{\sqrt{2}}{\pi}a\kappa^2$ when $a\kappa^2 \longrightarrow 0$.
If $a>0$ and $a\kappa^2$ is small enough we can conclude that $0 \not \in\sigma(P(a)). $\\
Let's look at the case $a<0$:\\
According to the min-max principle, the second eigenvalue of $P(a\kappa^2)$ is 
\begin{equation*}
\mu_1(a\kappa^2)=\inf_{u\in H,u\perp w_0}\int \left( |\nabla u|^2+(4\rho^2+\sigma^2)u^2+3a\kappa^2v_0^2u^2 \right)
-E_0(a\kappa^2),
\end{equation*}
and let $w_1$ realise the infimum.\\
We also have
\begin{equation*}
5=\inf_{u\in H,u\perp u_0}\int \left( |\nabla u|^2+(4\rho^2+\sigma^2)u^2 \right)=\inf_{u\in H,u\perp u_0}J(u,0),
\end{equation*}
realised for $u_1$, the second normalised Hermite function. Now, define $\tilde{u}=\alpha w_1+\beta w_0$ with $\alpha,\beta$ such that $\|\tilde{u}\|_{L^2}=\alpha^2+\beta^2=1$ and $\alpha\int w_1u_0+\beta\int w_1u_0=0$, then $\tilde{u}\in H$ and $\tilde{u}\perp u_0$. Notice that $|\alpha| \longrightarrow 1$ and $\beta \longrightarrow 0$ as $a\kappa^2 \longrightarrow 0$.\\
One has
$5=J(u_1,0)\leq J(\tilde{u},0),$
then  we obtain $5 \leq \mu_1(a\kappa^2)+\eps(a\kappa^2)$ with $\eps(a\kappa^2)  \longrightarrow 0$ as $a\kappa^2 \longrightarrow 0$, therefore $\mu_1(a\kappa^2) \geq 4$ for $a$ small enough, and $0 \not \in\sigma(P(a\kappa^2)). $\\
As a conclusion, for each choise of $E_1$, equation $(\ref{eqv_1})$ admits a solution $v_1\in L^2$ as the second right hand side $f$ is in $L^2$. However, if we choose $E_1$ so that $f \perp v_0$, we also have $\|v_1\|_{L^2}\leq C$ uniformly in $|a|\kappa^2 \leq c_0$, as the eigenvalue $E_0(a\kappa^2)$ is simple.\\
The estimations \eqref{estimv_j} are obtained as in the proof of Proposition \ref{propestimv_0}.\\
By the same argument we infer the existence of $v_2$ and $E_2$ which solve equation \eqref{eqv_2} 
and satisfy the estimates \eqref{estimv_j}.
\end{proof}

Take $\chi \in \mathcal{C}^{\infty}_0(\R)$ such that $\chi\geq 0$, supp$\chi \subset [\frac12,\frac32]$ and $\chi=1$ on $[\frac34,\frac54]$.\\
Set $v=\chi(\sqrt{h}\rho)(v_0+h^{\frac12}v_1+hv_2)$, $\tilde{v}(r,y,h)=v(\frac{r-k}{\sqrt{h}},\frac{y}{\sqrt{h}},h)$ and 
$\lambda=\frac{2k^2}{h}+E_0+h^{\frac12}E_1+hE_2$, and define 

\begin{equation}\label{uapp} 
u_{app}=\kappa h^{-\frac{1}{2}}\text{e}^{-i\lambda t}\text{e}^{i\frac{k^2}{h}\theta} \tilde{v}.
\end{equation}
Recall that, according to \eqref{dlE0}, $$E_0(a\kappa^2)=3+\frac{\sqrt{2}}{2\pi}a\kappa^2+o(a\kappa^2).$$
\begin{prop}\label{propuapp}  
The function $u_{app}$ defined by $(\ref{uapp})$ satisfies
\begin{equation}
ih\partial_{t}u_{app} +h^2\Delta u_{app}-|x|^2u_{app} =ah^2|u_{app}|^2u_{app} +R(h)
\end{equation}
with 
\begin{equation}\label{estR}
\|(|x|^2+1)R(h) \|_{L^2} \lesssim h^{\frac52}\quad \text{and} \quad\|\Delta R(h) \|_{L^2} \lesssim h^{\frac12}.
\end{equation}
\end{prop}

\begin{proof}
By construction, $w=v_0+h^{\frac12}v_1+hv_2$ satisfies $Eq(w)=h^{\frac52}R_1(h)$ where $Eq$ is 
defined by $(\ref{eqv})$, and according to Propositions \ref{propconv} and \ref{propuapp}

\begin{equation*}
|R_1(h)| \lesssim \left(\frac{1}{(k+\sqrt{h}\rho)^2}+|\rho|^3   \right)\e^{-c_1(|\rho|+|\sigma|)},
\end{equation*}
and

\begin{equation}\label{deltaR}
|\Delta R_1(h)| \lesssim \left(\frac{1}{(k+\sqrt{h}\rho)^4}+|\rho|^3   \right)\e^{-c_2(|\rho|+|\sigma|)}.
\end{equation}
Now,
\begin{eqnarray*}
Eq(v)&=&Eq(\chi(\sqrt{h}\rho)w)\\
      &=&\chi(\sqrt{h}\rho)Eq(w)-h \chi''(\sqrt{h}\rho)w-2h^{\frac12}\chi'(\sqrt{h}\rho)\partial_{\rho}w \\
      & & +a\chi(\chi^2-1)w^3\\
   &=& h^{\frac{5}{2}}\chi(\sqrt{h}\rho)R_1+R_2+R_3+R_4:=R(h).
\end{eqnarray*}
Set $I=[\frac12,\frac34] \cup[\frac54,\frac32] $ and observe that supp$\chi'\subset I$, supp$\chi''\subset I$,\\
 supp$\chi(\chi^2-1)\subset I$ and if $\sqrt{h}\rho \in I$ we have 
$$|w|, |\partial_{\rho}w|\lesssim \e^{-c/\sqrt{h}}\e^{-c|\sigma|}, $$
then it follows
\begin{equation}\label{deltaRi}
\|\Delta^j R_p \|_{L^2} \lesssim \e^{-c/\sqrt{h}},
\end{equation}
for all $0\leq j\leq 1$ and $2\leq p \leq 4.$
According to $(\ref{deltaR})$ we also have 
\begin{equation*}
\| \chi(\sqrt{h}\rho)R_1 \|_{L^2}^2 \lesssim \int (1+|\rho|^6)\e^{-2c_1(|\rho|+|\sigma|)} \leq C.
\end{equation*}
Therefore, coming back in variables $(r,y,\theta)$, $\|R(h)\|_{L^2} \lesssim h^{\frac52}$. Because of the fast decay  of $w$ we also have  $\|(r^2+y^2)R(h)\|_{L^2} \lesssim h^{\frac52}$, hence $\|(|x|^2+1)R(h) \|_{L^2} \lesssim h^{\frac52}$. \\Differentiating $u_{app}$ costs at most $h^{-1}$, then together with $(\ref{deltaR})$ and $(\ref{deltaRi})$ we obtain  $\|\Delta R(h) \|_{L^2} \lesssim h^{\frac12}.$
\end{proof}

\begin{prop}\label{propuapp2}
Let $|a|\kappa^2\leq c_0$ fixed, let $u_{app}$ be given by $(\ref{uapp})$ and let $u$ be solution of 
\begin{equation}\label{nls}
\left\{
\begin{aligned}
&i h\partial_t u+h^2\Delta u-|x|^2u   = ah^2|u|^{2} u,\\
&u(0,x)= u_{app}(0,x),
\end{aligned}
\right.
\end{equation}
then $\| (u-u_{app})(t_h)\|_{L^2} \longrightarrow 0$ with $t_h \ll \log(\frac1h)$,
 when $h \longrightarrow 0$.
\end{prop}

\begin{proof}
Denote by $w=u-u_{app}$ and by $f=ah^2g +R(h)$ with $g=|u_{app}+w|^{2}(u_{app}+w)-|u_{app}|^2u_{app}$, then
\begin{equation}\label{eqerr}
i h\partial_t w+h^2\Delta w-|x|^2w   = f.
\end{equation}
We define 
\begin{equation}\label{energie}
E(t)=\int \left(\frac12(|x|^4+1)|w|^2 +h^4 |\Delta w|^2\right).
\end{equation}

\begin{itemize}
\item Multiply $(\ref{eqerr})$ by $\frac12(|x|^4+1)\overline{w}$, integrate and take the imaginary part:
 \begin{equation}\label{eqen1}
\frac12 h\frac{\text{d}}{\text{d}t}\int \frac12(|x|^4+1)|w|^2
=\text{Im} \int \frac12(|x|^4+1)f\overline{w}+2h^2 \text{Im} \int|x|^2 \overline{w}x \nabla w,
\end{equation}

\item Multiply $\Delta (\ref{eqerr})$ by $h^4 \Delta\overline{w}$, integrate and take the imaginary part:
 \begin{equation}\label{eqen2}
\frac12h \frac{\text{d}}{\text{d}t}\int h^4|\Delta w|^2
=h^4\text{Im} \int \Delta f \Delta \overline{w}-2h^4 \text{Im} \int \Delta w  x \nabla \overline{w}.
\end{equation}
\end{itemize}
With an integration by parts, we can show that 
$$h^2\int |x|^2|\nabla w|^2\lesssim \int|x|^4|w|^2+h^4 \int |\Delta
w|^2, $$
therefore
 \begin{equation}\label{eqen3}
h^2\Big| \int  |x|^2 \overline{w}x \nabla w   \Big|\lesssim h\int
|x|^4|w|^2+h^3\int |x|^2|\nabla w|^2\lesssim hE,
\end{equation}
and
 \begin{equation}\label{eqen4}
h^4 \Big| \int \Delta w  x \nabla \overline{w} \Big|\lesssim h^5\int
|\Delta w|^2+h^3\int |x|^2|\nabla w|^2\lesssim hE.
\end{equation}
Then the inequalities \eqref{eqen1}-\eqref{eqen4}  yield
 \begin{equation}\label{eqE}
h \frac{\text{d}}{\text{d}t}E(t)\lesssim \text{Im} \int \left(\frac12(|x|^4+1)f\overline{w}+
h^2|x|^2\nabla f \nabla \overline{w}+h^4 \Delta f \Delta \overline{w}\right)+hE.
\end{equation}
Using the expression of $u_{app}$
\begin{eqnarray}
&\|u_{app}\|_{L^2} \lesssim 1, \quad  \|u_{app}\|_{L^{\infty}} \lesssim h^{-\frac12},\nonumber\\
&\|\nabla u_{app}\|_{L^2} \lesssim h^{-1},\quad  \|\nabla u_{app}\|_{L^{\infty}} \lesssim h^{-\frac32}\label{estuapp},
\end{eqnarray}
and by definition of $E$
\begin{equation}\label{est1w}
\|x\nabla w\|_{L^2} \lesssim h^{-1}E^{\frac12}, \quad \|\Delta w\|_{L^2} \lesssim h^{-2}E^{\frac12},
\end{equation}
and the Gagliardo-Nirenberg inequalities in dimension 3 yield
\begin{eqnarray}\label{est2w}
\| w\|_{L^4} \lesssim h^{-\frac34}E^{\frac12}, \quad \|\nabla w\|_{L^4} \lesssim h^{-\frac74}E^{\frac12},\nonumber\\
\| w\|_{L^{\infty}} \lesssim \| w\|_{L^2}^{\frac14}\|\Delta w\|_{L^2}^{\frac34} \lesssim h^{-\frac32}E^{\frac12}.
\end{eqnarray}

\begin{itemize}
\item
First, the estimates \eqref{estR} on $R(h)$ give

\begin{eqnarray}
&\left| \int \left(\frac12(|x|^4+1)R(h)\overline{w}+h^4 \Delta R(h) \Delta \overline{w}\right)\right|\nonumber\\
&\lesssim\|(|x|^2+1)R(h) \|_{L^2}E^{\frac12}
+ h^2 \|\Delta R(h)\|_{L^2}E^{\frac12} \nonumber\\
&\lesssim h^{\frac52}E^{\frac12}.\label{inegalite1}
\end{eqnarray}
\item
Then, as $g=|u_{app}+w|^{2}(u_{app}+w)-|u_{app}|^2u_{app}$, and 
according to \eqref{estuapp} and \eqref{est2w}

\begin{eqnarray}
\left|\text{Im} \int (|x|^4+1)g\overline{w}\right| & \lesssim & 
\int (|x|^4+1)\left( |u_{app}|^2|w|^2 + |u_{app}||w|^3   \right)\nonumber\\
& \lesssim & \|u_{app}\|_{L^{\infty}}(\|u_{app}\|_{L^{\infty}}+\|w\|_{L^{\infty}})E\nonumber\\
& \lesssim & h^{-1}E+h^{-2}E^{\frac32}.\label{inegalite2}
\end{eqnarray}

\item Compute
 
\begin{eqnarray}
|\Delta g| &\lesssim & |u_{app}|^2|\Delta w|+|u_{app}||\nabla u_{app}||\nabla w|+|\nabla u_{app}|^2|w|\nonumber\\
 & &+|u_{app}||\Delta u_{app}||w|  +|\Delta u_{app}||w|^2+|w|^2|\Delta w|+|w||\nabla w|^2,\nonumber
\end{eqnarray}
hence
\begin{eqnarray}
\|\Delta g \|_{L^2} &\lesssim & \|u_{app}\|_{L^{\infty}}^2\|\Delta w\|_{L^2}
+\|u_{app}\|_{L^{\infty}}\|\nabla u_{app}\|_{L^{\infty}}\|\nabla w\|_{L^2} \nonumber \\
& &+\|\nabla u_{app}\|_{L^{\infty}}^2\|w\|_{L^2}
+\|u_{app}\|_{L^{\infty}}\|\Delta u_{app}\|_{L^{\infty}}\|w\|_{L^2}\nonumber \\
& &  +\|\Delta u_{app}\|_{L^{\infty}}\|w\|_{L^4}^2
+\|w\|_{L^{\infty}}^2\|\Delta w\|_{L^2}+\|w\|_{L^2}\|\nabla w\|_{L^4}^2\nonumber \\
&\lesssim & h^{-3}E^{\frac12}+h^{-4}E+h^{-5}E^{\frac32},\nonumber
\end{eqnarray}
then
\begin{eqnarray}
h^4\left| \int \Delta g \Delta \overline{w}\right| & \lesssim & 
h^4\|\Delta g \|_{L^2}\|\Delta w \|_{L^2}\nonumber\\
 & \lesssim &  h^{-1}E+h^{-2}E^{\frac32}+h^{-3}E^2.\label{inegalite4}
\end{eqnarray}
\end{itemize}
Putting the  estimates  \eqref{inegalite1},  \eqref{inegalite2},  and  \eqref{inegalite4}
 together with $(\ref{eqE})$, we obtain
\begin{equation}\label{inegE}
h \frac{\text{d}}{\text{d}t}E(t) \lesssim h^{\frac52}E^{\frac12}+hE+E^{\frac32}+h^{-1}E^2.
\end{equation}
Set $F=E^{\frac12}$, then $F$ satisfies $F(0)=0$ and 
\begin{equation}\label{inegF}
h \frac{\text{d}}{\text{d}t}F(t) \lesssim h^{\frac52} +hF+F^2+h^{-1}F^3.
\end{equation}
As long as $h^{-1}F^3 \lesssim hF$, i.e. for times such that $F\lesssim h$, we can write
\begin{equation*}
 \frac{\text{d}}{\text{d}t}F(t) \lesssim h^{\frac32} +F.
\end{equation*}
Using Gronwall's inequality, $F \lesssim {h^{\frac32}}\e^{Ct}$. 
The non linear terms in $(\ref{inegF})$ can be removed with the continuity argument for times $t_h$ such that
$\e^{Ct_h} \lesssim h^{-\frac12}$, i.e. $t_h \ll \log(\frac1h)$ and one has $F(t_h)\longrightarrow 0$ when 
$h \longrightarrow 0$, hence the result.
\end{proof}

We are now able to prove Theorem $\ref{thmgeom}$ and Theorem $\ref{thmproj}$.

\section{Geometric instability}

Let $|a|\kappa^2\leq c_0$. Consider the function $u_{app}$ 
defined by \eqref{uapp} associate with $\kappa$ with $k=1$ ($k$ will be equal to $1$ in all this section).\\
\begin{equation*}
u_{app}=\kappa h^{-\frac12}\e^{-i\lambda t}\e^{i\frac{\theta}{h}}\tilde{v}.
\end{equation*}
Similarly, let the function  $u'_{app}$ defined by \eqref{uapp}
associate with $\kappa'=\kappa+h^{\frac12}$. Then there exists
$\lambda'\in \R$ and $\tilde{v}'\in L^2(\R^3)$ such that 
\begin{equation*}
u_{app}'=(\kappa+h^{\frac12}) h^{-\frac12}\e^{-i\lambda't}\e^{i\frac{\theta}{h}}\tilde{v}'.
\end{equation*}
define the functions $f,f'\in L^(\R^3)$ by 
\begin{equation}\label{deff1}
f=h^{-\frac12}\e^{i\frac{\theta}{h}}\tilde{v},\quad f'=h^{-\frac12}\e^{i\frac{\theta}{h}}\tilde{v}'.
\end{equation}
Notice that by construction, $\|f \|_{L^2}, \;\|f' \|_{L^2}\sim 1$.\\
We now need the following

\begin{lemm}\label{lemperturb}
The functions defined by \eqref{deff1}  satisfy
\begin{equation}\label{erreurfj}
\|f' -f\|_{L^2} \lesssim h^{\frac12}.
\end{equation}
\end{lemm}

\begin{proof}
To construct $f'$, we have to solve the system \eqref{eqv_0}-\eqref{eqv_2} with 
$\kappa'=\kappa+h^{\frac12}$. We reorganize this system by identifying the powers of $h$, and 
as equation \eqref{eqv_0} remains the same, we deduce \eqref{erreurfj}. 
\end{proof}

\begin{proof}[Proof of Theorem \ref{thmgeom} (i)]
Denote by $u$ (resp. $u'$) the solution of \eqref{nls} with initial condition $u_{app}(0)$ (resp. $u'_{app}(0)$).
We have 
\begin{eqnarray}\label{estinitiale}
\|(u'-u)(0) \|_{L^2}&=&\|(u_{app}'-u_{app}')(0) \|_{L^2}\nonumber\\
&\leq &\kappa \|f'-f  \|_{L^2}+\kappa h^{\frac12}   \|f'\|_{L^2} 
\lesssim \kappa h^{\frac12},
\end{eqnarray}
by Lemma \ref{lemperturb}. The triangle inequality gives 

\begin{eqnarray}\label{minoration}
\|(u_{app}'-u_{app})(t) \|_{L^2}  & \geq & 
\kappa \left|\e^{i(\lambda'-\lambda)t}  -1   \right| 
\|f'\|_{L^2}-\kappa \|f'-f  \|_{L^2}-  \kappa h^{\frac12}  \|f'\|_{L^2}\nonumber\\
 & \geq & \kappa \left|\e^{i(\lambda'-\lambda)t}  -1   \right|-C\kappa h^{\frac12} .
\end{eqnarray}
As $(\lambda'-\lambda)t\sim \frac{\sqrt{2}}{2\pi}a\left((\kappa+ h^{\frac12})^2-\kappa^2 \right)t 
\sim \frac{\sqrt{2}}{\pi}a\kappa t  h^{\frac12}$, with \eqref{minoration} we obtain, when $|a|\kappa t \gg 1$  
$$ \|(u_{app}'-u_{app})(t) \|_{L^2}\geq c |a| \kappa^2 t  h^{\frac12}, $$
hence, using \eqref{estinitiale}
$$\frac{ \|(u'-u)(t) \|_{L^2}}{ \|(u'-u)(0) \|_{L^2}} \gtrsim |a| \kappa t.$$
which was the claim.
\end{proof}

\begin{proof}[Proof of Theorem \ref{thmgeom} (ii)]
First notice that every parameter or function involved in this part depends on $h$ even 
though we do not write the subscripts. We define 
\begin{eqnarray}
u_{app}''&=&(\kappa+\eps_h) h^{-\frac12}\e^{-i\lambda''t}\e^{i\frac{\theta}{h}}\tilde{v}''\nonumber\\
         &:=&(\kappa+\eps_h) \e^{-i\lambda''t}f''. \label{deff3}  
\end{eqnarray}
with $\eps_h \longrightarrow 0$ when $h \longrightarrow 0$, and denote by $u''$ the solution of \eqref{nls} 
with initial condition $u_{app}''(0)$.
 Then
\begin{eqnarray}\label{estinitiale2}
\|(u''-u)(0) \|_{L^2}&=&\|(u_{app}''-u_{app})(0) \|_{L^2}\nonumber\\
&\leq &\kappa \|f''-f  \|_{L^2}+\kappa\eps_h   \|f''\|_{L^2} .
\end{eqnarray}
The right hand side of \eqref{estinitiale2} tends to $0$ with $h$ because $\|f''-f  \|_{L^2}\longrightarrow 0$ and $\|f'' \|_{L^2}\sim 1$. 
But when $h$ is small enough

\begin{eqnarray}\label{minoration2}
\|(u_{app}''-u_{app})(t) \|_{L^2}  & \geq & 
\kappa \left|\e^{i(\lambda''-\lambda)t}  -1   \right| 
\|f''\|_{L^2}-\kappa \|f''-f  \|_{L^2}-  \kappa\eps_h  \|f''\|_{L^2}\nonumber\\
 & \geq & \frac12\kappa \left|\e^{i(\lambda''-\lambda)t}  -1   \right|.
\end{eqnarray}
Now use  $(\lambda''-\lambda)t_h\sim \frac{\sqrt{2}}{2\pi}a\left((\kappa+\eps_h)^2-\kappa^2 \right)t_h 
\sim C_0a\kappa t_h \eps_h$.
Take $\eps_h=(C_0 \kappa a t_h)^{-1/2}$ which tends to $0$, then if $h\ll1$, $|\lambda''-\lambda|t_h \geq \pi$ and 
$$ \sup_{0\leq t \leq t_h} \|(u_{app}''-u_{app})(t) \|_{L^2}\geq \kappa.$$
Now, according to Proposition \ref{propuapp2}, which can be used as we assume $t\ll\log{\frac1h}$, we have for $h$ small enough
$$ \sup_{0\leq t \leq t_h} \|(u''-u)(t) \|_{L^2}\geq \kappa.$$
This last inequality together with \eqref{estinitiale2} proves the
second part of Theorem \ref{thmgeom}.
\end{proof}
\section{Projective instability}
We conserve the notations of the previous section, but here $f_j$ and $f'_j$ 
are constructed with $k=j$ in \eqref{defutilde}.\\
Define $U_{app}=\kappa \e^{-i\lambda_1t}f_1+\kappa \e^{-i\lambda_2t}f_2$ and 
$U'_{app}=(\kappa+\eps_h) \e^{-i\lambda'_1t}f'_1+\kappa \e^{-i\lambda_2t}f_2.$

\begin{lemm}\label{lemproj}
Let $V_{app}=U_{app}$ or  $V_{app}=U'_{app}$, and $v$ be  solution of 
                             
\begin{equation}\label{eqvapp}
\left\{
\begin{aligned}
&i h\partial_t v+h^2\Delta v-|x|^2v   = ah^2|v|^{2} v,\\
&v(0,x)= V_{app}(0,x),
\end{aligned}
\right.
\end{equation}
then $\| (v-V_{app})(t_h)\|_{L^2} \longrightarrow 0$ with $t_h \ll \log(\frac1h)$,
 when $h \longrightarrow 0$.
\end{lemm}

\begin{proof}
Write $V_{app}=v^1_{app}+v^2_{app}$ with $v^1_{app}=\kappa \e^{-i\lambda_1t}f_1$ or 
$v^1_{app}=(\kappa+\eps_h) \e^{-i\lambda'_1t}f'_1$ and $v^2_{app}=\kappa \e^{-i\lambda_2t}f_2$. 
As the supports of $v^1_{app}$ and $v^2_{app}$ are disjoint we have 
\begin{align*}
&i h\partial_t(v^1_{app}+v^2_{app}) 
+h^2\Delta(v^1_{app}+v^2_{app}) -|x|^2(v^1_{app}+v^2_{app}) \\
 &= ah^2\left(|v^1_{app}|^2v^1_{app}+ |v^2_{app}|^2v^2_{app}   \right)+R^1(h)+R^2(h)\\
 & =  ah^2|v^1_{app}+v^2_{app} |^2(v^1_{app}+v^2_{app})+R^1(h)+R^2(h),
\end{align*}
where for $j=1,2$, $R^j(h)$ is the error term given by Proposition \ref{propuapp} and therefore satisfies
 $\|(|x|^2+1)R^j(h) \|_{L^2} \lesssim h^{\frac52}$ and $\|\Delta R^j(h) \|_{L^2} \lesssim h^{\frac12}.$ 
We conclude with the help of Proposition \ref{propuapp2}.
\end{proof}

\begin{proof}[Proof of Theorem \ref{thmproj}]
Consider  the function $u$ (resp. $u'$ ) the solution of equation \eqref{eqvapp} with Cauchy data $U_{app}(0)$ 
(resp. $U'_{app}(0)$  ).\\
First notice that, for $t\geq 0$, $\| V_{app}(t)\|^2_{L^2} \sim 2\kappa^2$. Compute
\begin{equation}\label{produit}
U_{app}(t)\overline{U'_{app}}(t)=\kappa(\kappa+\eps_h)
f_1\overline{f'_1}\e^{i(\lambda'_1-\lambda_1)t}+\kappa^2|f_2|^2.
\end{equation}
Then for $t=0$ we have 
$$\int U_{app}\overline{U'_{app}}(0) \sim 2\kappa^2,$$
hence 
$$d_{\text{pr}}\left(u(0),u'(0)\right)=d_{\text{pr}}\left(U_{app}(0),U'_{app}(0)\right) \longrightarrow 0.$$
Let $t_h\ll \log{\frac1h}$, then as
$(\lambda'_1-\lambda_1)t_h\sim C_0a\kappa\eps_ht_h$, we now choose 
$$\eps_h=\frac{\pi}{C_0 a \kappa t_h},$$
then we have $(\lambda'_1-\lambda_1)t_h\longrightarrow \pi$, as
$h\longrightarrow 0$. Thus
$$\int U_{app}\overline{U'_{app}}(t_h) \longrightarrow 0,$$ and 
$$d_{\text{pr}}(U_{app}(t_h),U'_{app}(t_h)) \longrightarrow \arccos{(0)}=\frac{\pi}{2}. $$
Finally, from Lemma \ref{lemproj} we deduce 
 $$d_{\text{pr}}(u(t_h),U_{app}(t_h)), \;d_{\text{pr}}(u'(t_h),U'_{app}(t_h))\longrightarrow 0,$$
and therefore
\begin{eqnarray*}
d_{\text{pr}}(u(t_h),u'(t_h)) &\geq & d_{\text{pr}}(U_{app}(t_h),U'_{app}(t_h)) -d_{\text{pr}}(u(t_h),U_{app}(t_h))\\
                 &  &                                   -d_{\text{pr}}(u'(t_h),U'_{app}(t_h))\\
                  &\geq & \frac{\pi}{4},
\end{eqnarray*}
for $h \ll1;$ hence the result.
\end{proof}



\nocite{CCT1}
\nocite{CCT2}

\backmatter

\backmatter

\end{document}